\documentclass[a4paper,20pt]{article}

\usepackage{amsthm}
\usepackage{caption}
\usepackage{subfigure}
\usepackage{amssymb}
\usepackage{graphicx,subfig}
\usepackage{color}
\usepackage{amsmath}
\usepackage{xcolor}
\usepackage{colortbl,booktabs}
\usepackage{ebezier}
\usepackage{geometry}

\newtheorem{theorem}{Theorem}
\newtheorem{lemma}{Lemma}
\newtheorem{proposition}{Proposition}
\newtheorem{definition}{Definition}

\newtheorem{remark}{Remark}

\begin{document}

\title{Minimum leader selection for Structural Controllability of Undirected Graphs with Leader-follower Framework}
\author{Li Dai\\
College of Liberal Arts and Sciences,\\
National University of Defense Technology, Changsha 410072, China\\
 daili@nudt.edu.cn}

\maketitle

%\LARGE
\textbf{abstract}
The optimization problem of the minimum set of leaders for the controllability of undirected graphs are addressed.  It is difficult to find not only its optimal solution but also its approximate algorithm.
We propose a new concept, namely minimal perfect critical set (MPCS), to obtain an optimal solution. Some properties are presented, and
on the basis of these theorems, the problem of the minimum set of leaders of two typical self-similar bipartite networks, namely deterministic scale-free networks (DSFN) and Cayley trees, is solved completely.

\textbf{keyword}
Controllability, undirected graph; minimum leader selection; minimal perfect critical set; tree graph

\section{Introduction}

In recent years, the controllability of undirected networks has attracted considerable attention.
In contrast to the various insights and results on the controllability of directed networks, there remain many basic problems that need to be investigated with regard to the controllability of undirected networks. For example, to ensure the controllability of an undirected network, it is necessary to determine how to select the leader (leaders), what is the minimum number of leaders, how to get the minimum  set of leaders, and how large is the minimum set of leaders. These problems are interesting because of their potential application to various fields such as complex networks, multi-agent systems, and biology.
 \subsection{Literature Review}
Our goal is to address the problem of how to select leaders and how to obtain the minimum number and minimum set of leaders of Laplacian-based leader-follower undirected systems. Therefore, we review only the related studies. The networks mentioned throughout this paper are undirected networks, unless stated otherwise.

Research on the controllability of undirected network began with the work of Tanner \cite{Tanner}, where both a leader-follower framework and a Laplacian matrix were considered to understand how the network structure and leader's location affect the controllability of a network. Many studies have been conducted since then. The methods used in this field mainly include spectral analysis, equitable partitions, zero forcing sets, balancing sets (see \cite{ShimaAndMohammad}), and constraint matching.

The necessary or sufficient algebraic conditions of the leader set in view of eigenvectors
were presented for the controllability of a general network in \cite{Meng}-\cite{Zhijian}. Owing to the existence of multiple eigenvalues and the need to examine eigenvectors with zero components, it is computationally difficult to directly use the algebraic conditions to find the set of leaders. Therefore, network controllability based on the leader's graphical structure has attracted widespread attention.

In \cite{Rahmani&Mesbahi}, the authors applied the concept of equitable partitions to the study of controllability and pointed out that some vertex sets with a symmetric character should not be selected as leaders. The necessary graph-theoretic conditions for controllability involving
equitable partitions or almost equitable partitions were subsequently reported; see, for instance, \cite{Martini}-\cite{Cesar}.

The necessary/sufficient conditions for controllability based on (generalized) zero forcing sets have been established in the literature \cite{AIM}-\cite{Yashashwi}. In these studies, a one-to-one correspondence was established between the set of leaders rendering the
network controllable and the zero forcing sets.

Some studies, such as \cite{Trefois},\cite{Olesky}-\cite{Airlie}, have focused on constraint matching to investigate  controllability. A single leader and multiple leaders have been considered in \cite{Olesky} and \cite{Airlie}, respectively. The relationship between zero forcing sets and constraint matching was discussed in \cite{Trefois}.

Special topologies have been considered in some studies, such as trees \cite{JiAndLinAndLee}-\cite{JiAndLinAndYu}, grids \cite{Parlangeli}, and circulants \cite{Nabi-abdolyousefi}.

The above-mentioned studies constitute the profound literature on the theoretical research of controllability. Some other studies are devoted to algorithm design.
 As pointed out in several studies, finding a minimal zero forcing set (see \cite{Trefois} ) or finding a minimum leader set (see\cite{Airlie}) for a general network is an NP-hard problem, as is finding an approximate solution (see \cite{Yashashwi}). Moreover, there are several different concepts of controllability in the field of network controllability, such as controllability as defined by Kalman in \cite{Kalman}, structural controllability (defined by Lin in \cite{Lin}), and strong structural controllability (defined by Mayeda and Yamada in \cite{Mayeda}). The matrix can be defined as the adjacency matrix or Laplacian matrix with variable or fixed parameters. Thus, research on the controllability of networks is extensive.
 
Different from the existing studies, this paper proposes the concept of a minimal perfect critical set (MPCS) in section \ref{section2}. In addition, it presents the necessary and sufficient conditions for $S$ to be a minimal perfect critical 2 set and proves an interesting result that a minimal perfect critical 3 set does not exist.
Section \ref{section3} constitutes the main part of this paper. It presents some special MPCS of tree graphs.
Section \ref{section4} presents three examples to illustrate  how to use the theories proposed to find MPCS . Finally, Section \ref{section5} concludes the paper.

\subsection{Notations and Preliminary Results}
Let $G=(V,E)$ be an undirected and unweighted simple graph, where $V=\{v_1,v_2,\cdots,v_n\}$ is a vertex set and $E=\{v_iv_j|v_i\,\, and \,\,v_j\in V\}$ is an edge set, where an edge $v_iv_j$ is an unordered pair of distinct vertices in $V$. If $v_iv_j\in E$, then $v_i$ and $v_j$ are said to be \textit{adjacent} or \textit{neighbors}.
$N_S(v_i)=\{v_j\in S| v_iv_j\in E(G)\}$ represents the neighboring set in $S$ of $v_i$, where $S\subset V$.
The cardinality of $S$ is denoted by $|S|$. $d_G(v)$ is the degree of vertex $v$ in $G$ and $d_G(v)=|N_V(v)|$. If $d_G(v)=1$, then $v$ is called a pendent vertex.
$G[S]$ is the induced subgraph, whose vertex set is $S$ and edge set is $\{v_iv_j\in E(G)|v_i,v_j\in S\}$.
The \textit{valency matrix} $\Delta(G)$ of graph $G$ is a diagonal matrix with rows and columns indexed by $V$, in which the $(i,i)$-entry is the degree of vertex $v_i$, e.g., $|N_V(v_i)|$. Any undirected simple graph can be represented by its \textit{adjacency matrix}, $D(G)$, which is a symmetric matrix with 0-1 elements. The element in position $(i,j)$ in $D(G)$ is 1 if vertices $v_i$ and $v_j$ are adjacent and 0 otherwise. The symmetric matrix defined as
\[\mathbf{L}(G)=\mathbf{\Delta}(G)-\mathbf{D}(G)\]
is the \textit{Laplacian } of $G$.
Throughout this paper, it is assumed without loss of generality that $F$ denotes a follower set; its vertices  play the role of followers and the vertices in $\overline{F}$  are leaders(driver nodes), where $\overline{F}=V \backslash F$ denotes the complement set of $F$.
Let $\mathbf{y}$ be a vector and $\mathbf{y}_S$ denote the vector obtained from $\mathbf{y}$ after deleting the elements in $\overline{S}$.
Let $\mathbf{L}_{S\rightarrow T}$ denote the matrix obtained from $\mathbf{L}$ after deleting the rows in $\overline{S}$ and columns in $\overline{T}$.
If the followers' dynamics is (see (4) in \cite{Tanner})
\[\dot{\mathbf{x}}=\mathbf{Ax}+\mathbf{Bu},\]
where $\mathbf{x}$ is the   vector of all  $x_i$ corresponding to follower  $v_i\in F$ and $\mathbf{u}$ is the external control
input vector that is injected to only leaders,
then, the system is controllable if and only if the
$N \times NM$ controllability matrix
\[\mathbf{C}=[\mathbf{B},\mathbf{AB},\mathbf{A^2B},\cdots,
\mathbf{A^{N-1}B}]\]
has full row rank,
where $\mathbf{B}=\mathbf{L}_{F\rightarrow \overline{F}}$ and $\mathbf{A}=\mathbf{L}_{F\rightarrow F}$.
This   mathematical condition for controllability is well known as Kalman's controllability rank condition \cite{Kalman,Yangyu,Brockett}.
 In terms of eigenvalues and eigenvectors of Laplacian submatrices, a necessary and sufficient algebraic condition on controllability was presented in \cite{Meng} and \cite{Zhijian}(see Proposition 2,1(2) in \cite{Meng} and Proposition 1 in \cite{Zhijian} ).
\begin{proposition}\label{proposition2}
\cite{Meng,Zhijian}
The undirected graph $G$ is controllable under the leader set $\overline{F} $ if and only if $\mathbf{y}|_{\overline{F}}\neq \mathbf{0}$ \, (\,\,$\forall\,\,\mathbf{y}$ is an eigenvector of $\mathbf{L}$).
\end{proposition}
According to Proposition \ref{proposition2}, for any $S\subset V$ and $S\neq \emptyset$,
if there exists an eigenvector $\textbf{y}$ of Laplacian matrix $\textbf{L}$ such that $\mathbf{y}_{S}=\mathbf{0}$, then $S$ cannot be used as a leader set.
\section{Minimal Perfect Critical Set}\label{section2}
For a vertex set $S$, from Proposition \ref{proposition2}, if there exist an eigenvector $\textbf{y}$ such that $\mathbf{y}_{\overline{S}}=\mathbf{0}$, then $\overline{S}$ should not be used as a leaders set. In other words, there must be a vertex $v$ in $S$ to be selected as a member of the leaders set. Therefore, to locate the leaders   of graph $G$, three new concepts are proposed.
\begin{definition}\label{definition1}
(\textbf{critical set})\,\,Let $S$ be a nonempty subset of $V$. If there exists an eigenvector $\textbf{y}$ such that $\textbf{y}_{\overline{S}}=\textbf{0}$, then $S$ is called a critical set(CS) and $\textbf{y}$ is a inducing eigenvector.  $S$ is called a critical $k$  set if $|S|=k$.
\end{definition}
\begin{definition}\label{definition2}
(\textbf{perfect critical set})\,\,
Let $S$ be a critical set. If there exists an eigenvector $\textbf{y}$ satisfying $\textbf{y}_{v_i}\ne 0 (\forall v_i\in S)$, then $S$ is called a perfect critical set (PCS). $S$ is called a perfect critical $k$  set if $|S|=k$.
\end{definition}
\begin{definition}\label{definition 3}
(\textbf{minimal perfect critical set})
A perfect critical set is called a minimal perfect critical set (MPCS) if any proper subset of it is no longer a PCS. $S$ is called a minimal perfect critical $k$  set if $|S|=k$.
\end{definition}
From the definitions stated above, $V$ is a trivially CS and a PCS induced by the eigenvector $\textbf{1}_n$. $V$ is an MPCS if and only if $G$ is controllable under any single vertex selected as a leader; e.g., $G$ is omnicontrollable.
\begin{figure}[h]
\centering
\scriptsize
\setlength{\unitlength}{1mm}
\begin{picture}(40,32)
\put(0,5){\circle{1.5}}
\put(30,5){\circle{1.5}}
\put(0,20){\circle{1.5}}
\put(0,34.5){\circle{1.5}}
\put(15,20){\circle{1.5}}
\put(43,18.5){\circle{1.5}}
\put(27.8,33.2){\circle{1.5}}
\put(-5,4){$v_1$}
\put(-5,19){$v_2$}
\put(-5,34){$v_3$}
\put(13,23){$v_4$}
\put(32,4){$v_7$}
\put(42.5,20){$v_6$}
\put(28,35){$v_5$}
\put(0.6,5.5){\line(1,1){14}}
\put(15.4,20.6){\line(1,1){12}}%v4-v5
\put(0.6,20){\line(1,0){13.8}}%v2-v4
\put(14.5,20){\line(-1,1){14}}%v4--v3
\put(42.5,19){\line(-1,1){14}}%v6--v5
\put(29.5,5.5){\line(-1,1){14}}%v7--v4
\put(30.5,5.5){\line(1,1){12.3}}%v7--v4
\end{picture}
\caption{Graph $G$ with 4 distinct MPCS}
\label{figureForDifferentMPCS}
\end{figure}
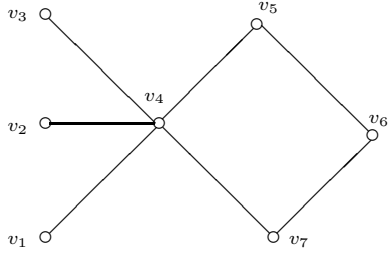
For example, see Fig.\ref{figureForDifferentMPCS}. $\{v_1,v_5\}$ is not a CS; $\{v_1,v_2,v_3,v_4\}$ is a CS but not a PCS. $S=\{v_1,v_2,v_3\}$ is a PCS but not an MPCS. $S_1=\{v_1,v_2\}$, $S_2=\{v_1,v_3\}$, $S_3=\{v_2,v_3\}$, and $S_4=\{v_5,v_7\}$ are all MPCS.
\begin{remark}\label{remarkForPropositionStatedWithMPVCS}
 Proposition \ref{proposition2} can be restated as follows: The undirected graph $G$ is controllable under the leader set $\overline{F} $ if and only if, for each MPCS $S$ , $S \bigcap \overline{F}\ne \emptyset$.
\end{remark}
From  Remark \ref{remarkForPropositionStatedWithMPVCS} , there is a close relationship between MPCS and the minimum leader set.
In other words, when we find all MPCS of $G$, we find the minimum leader set and hence the minimum number of leader vertices. For example, by Remark \ref{remarkForPropositionStatedWithMPVCS} and the 4 MPCS stated above, graph $G$ in Fig.\ref{figureForDifferentMPCS} is not controllable under any single leader because $\bigcap_{i=1}^{4}S_i=\emptyset$ , or under any two vertices. Therefore, the minimum leader set is $\{v_i,v_j,v_k\}$, where $v_i,v_j$ comes from $\{v_1,v_2,v_3\}$ and $v_k$ comes from $\{v_5,v_7\}$. Therefore, the minimum number of leaders is 3.

Moreover, many MPCS have typical graphical characteristics. For example, all of the 4 MPCS of $G$ in Fig.\ref{figureForDifferentMPCS} have the graphical structure stated in  Theorem \ref{theorem3}.

\subsection{Sufficient conditions for Critical Set}\label{subsection2.2}
For an undirected connected graph, Laplacian matrix $\textbf{L}$ is  symmetric and all the eigenvectors corresponding to different eigenvalues are orthogonal to each other. Hence, by knowing that $\textbf{1}_n$ is an eigenvector of $\textbf{L}$, it is immediate that all the other eigenvectors of $\textbf{L}$ are orthogonal to $\bold{1}_n$, i.e., for all eigenvectors $ \bold{y}$,
\begin{equation}\label{eq1}
\textbf{1}_{n}^T\textbf{y}=\sum_{i=1}^{n}y_i=0.
\end{equation}
For all nontrivial graphs, if $S$ is a CS, then
\begin{equation}\label{equationFor|S|>=2}
|S|\geq 2.
\end{equation}
 In fact,
suppose that $|S|=1$; without loss of generality, $S=\{v_1\}$.
Let $\textbf{y}=(y_1,y_2,\cdots,y_n)^T$ be the inducing eigenvector  associated with eigenvalue $\lambda$. Then, $\textbf{Ly}=\lambda \textbf{y}$ and $\textbf{y}_{\overline{S}}=\textbf{0}$.
By $\textbf{y}_{\overline{S}}=\textbf{0}$ and (\ref{eq1}), $\textbf{y}_{S}=\textbf{0}$, e.g., $\textbf{y}=\textbf{0}$.
Further, since any subset $S$ with $|S|=1$ is not a CS, by Remark \ref{remarkForPropositionStatedWithMPVCS}, $G$ is controllable with the leader set $\overline{F}$ when $|F|=1$.
Now, we investigate the properties of a CS. First, a sufficient condition for $S$ to be a CS is provided in the following Proposition \ref{proposition3}, which describes a special case of the symmetry-based uncontrollability results.
\begin{proposition}\label{proposition3}
Let $G$ be an undirected connected graph of order $n$, $S\subset V$, and $|S|\geq 2$. If for any $v\in \overline{S}$, either $N_S(v)=\emptyset$ or $N_S(v)=S$, then $S$ is a critical set.
\end{proposition}
\textbf{Proof}\quad Let $|\{v\in \overline{S}|N_S(v)=S\}|=m$. Then,
$\textbf{L}_{S\rightarrow S}-m\textbf{I}_{|S|}$
is the Laplacian of subgraph $G[S]$, where $\textbf{I}_{|S|}$ denotes the $|S|$-dimensional identity matrix. Considering (\ref{eq1}), there exists an eigenvector $\textbf{y}_0$ of the Laplacian $\textbf{L}_{S\rightarrow S}-m\textbf{I}_{|S|}$ such that $\textbf{1}_{|S|}^T\textbf{y}_0=0$.

Set vector $\textbf{y}$ as $\textbf{y}_S=\textbf{y}_0$ and $\textbf{y}_{\overline{S}}=\textbf{0}$. It can be seen that
 \begin{equation}\label{eq2}
 \textbf{Ly}=\left[ \begin{array}{ll}
  \textbf{L}_{S\rightarrow S}& \textbf{L}_{S\rightarrow \overline{S}} \\
  \textbf{L}_{\overline{S}\rightarrow S}& \textbf{L}_{\overline{S}\rightarrow \overline{S}}
  \end{array}\right]
  \left[\begin{array}{c} \textbf{y}_0\\ \textbf{0 }\end{array}\right]=\left[\begin{array}{c}  \lambda \textbf{y}_0\\  \textbf{L}_{\overline{S}\rightarrow S}\textbf{y}_0 \end{array}\right].
 \end{equation}
Noting that the rows in matrix $\textbf{L}_{\overline{S}\rightarrow S}$ are either  ones or zeros, the conclusion is proved by $\textbf{1}_{|S|}^T\textbf{y}_0=0$(see (\ref{eq1})).   \qed
\subsection{Minimal Perfect Critical 2  and  3  Set }
Based on the above-mentioned properties, the critical $k$  set with $k\leq 3$ can be determined directly from the graphical characterization.
This is achieved via a detailed analysis of the inducing eigenvector.
\begin{lemma}\label{lemma1}
Let $G$ be an undirected connected graph and $S$ be a perfect critical $k$ vertices set. Then, for any $v\in \overline{S}$, $|N_S(v)|\ne 1$ and $|N_S(v)|\ne k-1$.
\end{lemma}
\textbf{Proof}\quad Let $S=\{v_1,v_2,\cdots,v_k\}$ be a perfect critical vertex set and $\textbf{y}=(y_1,y_2,\cdots,y_k,0,0,\cdots,0)^T$ be the inducing eigenvector. $y_i\ne 0(\forall 1\leq i \leq k)$ since $S$ is a perfect critical vertex set.
$\forall v\in \overline{S}$, suppose that $|N_S(v)|=1$. Without loss of generality, say, $vv_1\in E$ and $vv_i\notin E(G)(\forall i\neq 1)$. Then, $\textbf{L}_{\{v\}\rightarrow V}\textbf{y}=y_1\neq 0$. On the other hand, $\textbf{y}|_{\overline{S}}=\textbf{0}$ and $v \in \overline{S}$. Hence,
$\textbf{L}_{\{v\}\rightarrow V}\textbf{y}=0$. This is a contradiction.
Together with (\ref{eq1}), $|N_S(v)|\ne k-1$ can be proved similarly.\qed

By (\ref{equationFor|S|>=2}), a critical 2 set is also a minimal perfect critical 2 set. The following Theorem \ref{theorem3} will follow from Lemma \ref{lemma1} and Proposition \ref{proposition3}.

\begin{theorem}\label{theorem3}
Let $G$ be an undirected connected graph, $S\subset V$, and $|S|=2$. Then, $S$ is a minimal  perfect critical 2  set if and only if $\forall v\in \overline{S}$, either $N_S(v)=\emptyset$ or $N_S(v)=S$.\qed
\end{theorem}

For example, see graph $G$ in Fig.\ref{figureForDifferentMPCS}.
All its 4 MPCS  can be recognized by the graphical characterization stated in Theorem \ref{theorem3}.

 Next, we will prove that a minimal perfect critical 3 set does not exist in the following theorem.

\begin{figure}[h]
\centering
\setlength{\unitlength}{1mm}
\begin{picture}(80,12)
\put(0,5){\circle*{1.5}}
\put(10,5){\circle{1.5}}
\put(25,5){\circle{1.5}}
\put(35,5){\circle*{1.5}}
\put(50,5){\circle*{1.5}}
\put(60,5){\circle{1.5}}
\put(75,5){\circle*{1.5}}
\put(85.5,5){\circle{1.5}}

\put(25.5,5.5){\line(1,1){4.1}}
\put(5.2,10.5){\circle{1.5}}
\put(30,10.2){\circle{1.5}}
\put(50.5,5.5){\line(1,1){4.1}}
\put(55,10.2){\circle{1.5}}
\put(50.5,5){\line(1,0){9}}
\put(75.5,5.5){\line(1,1){4.1}}
\put(80.2,10.2){\circle{1.5}}
\put(80.7,9.6){\line(1,-1){4.1}}
\put(75.7,5){\line(1,0){9}}

\end{picture}
\vspace*{-0.4cm}
\caption{All possible topology structures of $G[S]$ with $|S|=3$}
\label{AllPossibleTopologyStructures}
\end{figure}
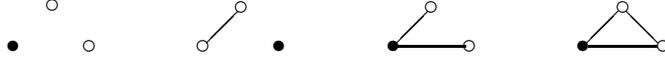

\begin{theorem}\label{theorem4}
Let $G$ be an undirected connected graph, $S\subset V$, and $|S|=3$. Then, $S$ is not a minimal perfect critical set.
\end{theorem}
\textbf{Proof}\quad Suppose that $S$ is a minimal perfect critical set. Consider the subgraph $G[S]$. All  4  possible topology structures of $G[S]$ are shown in Fig. \ref{AllPossibleTopologyStructures}.

%%%%%%%%%%%  ? 3      %%%%%%%%%%%
%%%%%%%%%%%% ? 3  ? ?    %%%%%%%%%

For each topology of $G[S]$ in Fig.\ref{AllPossibleTopologyStructures}, let $T$ be the vertex set of white vertices and $u$ be the black vertex. One can have either $N_T(u)=\emptyset$ or $N_T(u)=T$.

 By Lemma \ref{lemma1}, $\forall v\in\overline{S}$, either $|N_S(v)|=0$ or $|N_S(v)|=3$.

Noting that $T\subset S$ and $\overline{T}=\{u\}\bigcup\overline{S}$, by Proposition \ref{proposition3}, $T$ is a critical set. This is in contradiction with $S$ as a MPCS.\qed
\begin{figure}[h]
\centering
\setlength{\unitlength}{1mm}
\begin{picture}(120,30)
\put(20,5){\circle{1.5}}
\put(20,5.5){\line(0,1){15}}
\put(20,21){\circle*{1.5}}
\put(19.5,5.5){\line(-1,1){15}}
\put(4.5,21){\circle*{1.5}}
\put(20.5,5.5){\line(1,1){15}}
\put(35.5,21){\circle*{1.5}}
\put(35.5,5){\circle{1.5}}
\put(35.5,5.5){\line(0,1){15}}
\put(36,5.5){\line(1,1){15}}
\put(51,21){\circle*{1.5}}
\put(20.5,5.5){\line(2,1){30}}
\put(1,12){\dashbox{2}(55,18)}
\put(35.5,20.5){\oval(31,15)[t]}
\put(20,20.5){\oval(31,8)[t]}
\put(4,26){\scriptsize $S_1$}
\put(1,-2){\scriptsize  $(a)\,\, S_1$ is minimal perfect critical set}

\put(90,5){\circle{1.5}}
\put(90,5.5){\line(0,1){15}}
\put(90,21){\circle*{1.5}}
\put(89.5,5.5){\line(-1,1){15}}
\put(74.5,21){\circle*{1.5}}
\put(90.5,5.5){\line(1,1){15}}
\put(105.5,21){\circle*{1.5}}
\put(105.5,5){\circle{1.5}}
\put(105.5,5.5){\line(0,1){15}}
\put(106,5.5){\line(1,1){15}}
\put(121,21){\circle*{1.5}}
\put(90.5,5.5){\line(2,1){30}}

\put(71,12){\dashbox{2}(55,18)}
\put(105.5,21.5){\oval(31,9)[t]}
\put(74,25){\scriptsize$S_2$}
\put(76,-2){\scriptsize $(b)\,\, S_2$ is not a critical set}
\end{picture}
\vspace{0.1cm}
\caption{Topology structure  effect on MPCS when $|S|\geq 4$}
\label{PerfectCriticalVertexSetIsCloselyRelatedTo}
\end{figure}
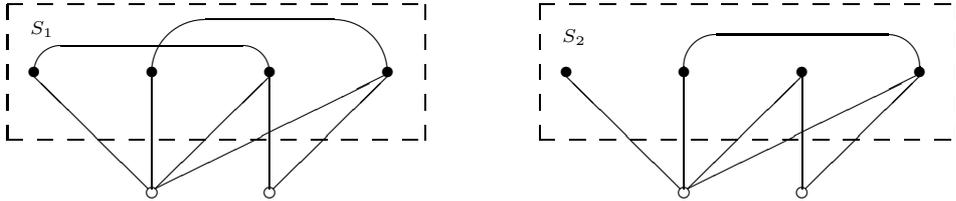

Although a minimal perfect critical 3 set does not exist, a minimal perfect critical 4 set does exist; see Fig.\ref{PerfectCriticalVertexSetIsCloselyRelatedTo}(a).
The topology structure of $G[S]$ will affect whether $S$ is an MPCS; see Fig.\ref{PerfectCriticalVertexSetIsCloselyRelatedTo}.
\section{Minimal Perfect Critical Set of Trees}\label{section3}
In this section, we first introduce a method for the simplification of $G$ and then prove some properties of MPCS of trees.
\subsection{Simplification Graph $\tilde{G}$}
Let $S$ be a proper subset  of $V(G)$. Let $\tilde{G}$ denote a graph constructed in the following steps:

Initialize $\tilde{G}=G$.

\textcircled{1}If $v\in\bar{S}$ and  $N_G(v)\cap S=S$, then set $\tilde{G}:=\tilde{G}-v$.

\textcircled{2}If $v\in \bar{S}$ and  $N_G(v)\cap S=\emptyset$, then set
$\tilde{G}:=\tilde{G}-v$.

 $S$ is a CS of a connected graph $G$ if and only if $S$ is a CS of $\tilde{G}$. Further, $\lambda$ is the corresponding eigenvalue of $\textbf{L}(G)$ if and only if $\lambda-1$,% ?????1 ???S????????1%
$\lambda$ is an eigenvalue of $\textbf{L}(\tilde{G})$ in case \textcircled{1}, \textcircled{2}, respectively.

\subsection{Necessities of MPCS of Trees}
\begin{proposition}\label{proposition4}
If $S$ is a CS of tree graph $G$, then there are at least 2 pendent vertices in $S$.
\end{proposition}
\textbf{Proof} Let $v_0\in \bar{S}$ be a pendent vertex in tree $G$ and $u_0$ be the vertex adjacent to $v_0$. By Lemma \ref{lemma1}, we  have $u_0\in \bar{S}$.
Hence, $S$ is still a PCS of tree graph $G-v_0$. Set $\bar{S}_1:=\bar{S}-v_0$, $G_1:=G-v_0$. It is useful to note that $d_G(v)=d_{G_1}(v)(\forall v\in S)$. Repeat these steps to obtain a tree graph $G'$ with $d_{G'}(v)>1(\forall v\in \bar{S})$. It is well known that there are at least two pendent vertices in a tree graph. Considering  the tree graph $G'$, we know that the pendent vertices are all in $S$. This completes the proof.   \qed

By Remark \ref{remarkForPropositionStatedWithMPVCS} and Proposition \ref{proposition4}, we have the following theorem:

\begin{theorem}\label{theoremForLeafivesLeaders}
Let $S_0=\{v\in V(G)| d_G(v)=1\}$ and $F$ be a follower set. If $|F\cap S_0|\leq 1$, then tree graph $G$ is controllable under leader set $\bar{F}$.\qed
\end{theorem}

Theorem \ref{theoremForLeafivesLeaders} tells us that when $|S_0|-1$ vertices in $S_0$ are selected as leaders, the graph will be controllable. Therefore, the minimum number of leaders is no more than $|S_0|-1$.

Next, we will study the property that the MPCS of trees should have.
\begin{proposition}\label{propositionForG[S]Disconnected}
If $S$ is a nontrivial PCS of tree graph $G$, then $G[S]$ is a disconnected subgraph.
\end{proposition}
\textbf{Proof}\quad Suppose that $G[S]$ is a connected subgraph. Arbitrarily take a vertex $v$ in $\bar{S}$. By Lemma \ref{lemma1}, there exist at least two vertices in $S$, say $u_1,u_2$ such that $vu_i\in E(G)(i=1,2)$. Since $G[S]$ is connected, there exists a circle $C$ with $u_1,u_2,v\in V(C)$ in $G$, and this is a contradiction to $G$ being a tree.   \qed

\begin{figure}[h]
\centering
\setlength{\unitlength}{1mm}
\begin{picture}(100,60)
\color{green}
\put(0,38){\dashbox{2}(100,20)}
\put(0,0){\dashbox{2}(100,20)}
\color{black}
\put(2,46){\scriptsize $S$}
\put(2,10){\scriptsize $\bar{S}$}
\put(8,45){\dashbox{1}(8,6)}
\put(10,48){\circle*{1.5}}
\put(25,45){\dashbox{1}(8,6)}
\put(27,48){\circle*{1.5}}
\put(30,48){\circle*{1.5}}
\put(18,48){\line(1,0){1}}
\put(20,48){\line(1,0){1}}
\put(22,48){\line(1,0){1}}
\put(38,45){\dashbox{1}(8,6)}
\put(41,48){\circle*{1.5}}
\put(54,45){\dashbox{1}(9,6)}
\put(56,48){\circle*{1.5}}
\put(61,48){\circle*{1.5}}
\put(67,45){\dashbox{1}(9,6)}
\put(69,48){\circle*{1.5}}
\put(74,48){\circle*{1.5}}
\put(78,48){\line(1,0){1}}
\put(80,48){\line(1,0){1}}
\put(82,48){\line(1,0){1}}
\put(85,45){\dashbox{1}(8,6)}
\put(88,48){\circle*{1.5}}
\color{red}
\put(50,54){\line(0,1){8}}
\put(50,46){\line(0,1){5}}
\put(50,38){\line(0,1){5}}
\put(50,30){\line(0,1){5}}
\put(50,22){\line(0,1){5}}
\put(50,14){\line(0,1){5}}
\put(50,6){\line(0,1){5}}
\put(50,-2){\line(0,1){5}}
\put(50,-6){\line(0,1){3}}
\color{black}
\color{blue}
\put(6,41){\dashbox{1.5}(42,15)}
\put(53,41){\dashbox{1.5}(42,15)}
\color{black}

\put(36,53){\scriptsize{$S_u$}}
\put(56,53){\scriptsize{$S_1$}}

\put(8,5){\dashbox{1}(8,6)}
\put(10,8){\circle*{1.5}}
\put(35,3){\dashbox{1}(24,11)}
\put(45,4.5){\scriptsize{$T$}}
\put(37,8){\circle*{1.5}}
\put(56,8){\circle*{1.5}}
\put(48,8){\circle*{1.5}}
\put(47,10){\scriptsize{$w$}}
\put(54,10){\scriptsize{$v$}}
\put(36,10){\scriptsize{$u$}}
\put(48,8){\line(1,0){8}}
\put(37,8){\line(1,0){3}}
\put(51,7.25){\scriptsize{$\times$}}
\put(41,8){\line(1,0){1}}
\put(43,8){\line(1,0){1}}
\put(45,8){\line(1,0){1}}
\put(20,8){\line(1,0){1}}
\put(25,8){\line(1,0){1}}
\put(30,8){\line(1,0){1}}
\put(68,8){\line(1,0){1}}
\put(73,8){\line(1,0){1}}
\put(78,8){\line(1,0){1}}
\put(85,5){\dashbox{1}(8,6)}
\put(88,8){\circle*{1.5}}
\cbezier[500](37,8)(39,7)(40,30)(41,47.5)
\cbezier[500](37,8)(35,7)(28,30)(27.3,47.5)
\cbezier[500](10,8)(12,9)(25,30)(26.7,47.5)
\cbezier[500](56,8)(57,9)(60,40)(69,47.5)
\put(10,8){\line(0,1){40}}
\put(56,8){\line(0,1){40}}
\put(88,8){\line(-2,3){27}}
\put(88,8){\line(0,1){40}}
\end{picture}
\vspace{0.5cm}
\caption{Component $T$ of $\tilde{G}[\bar{S}]$ with at least two vertices}
\label{ProofForBarSIsolatedSet}
\end{figure}
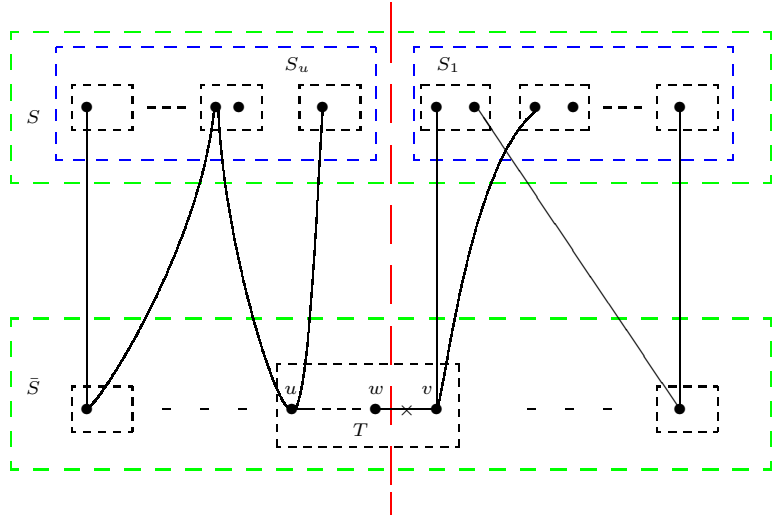

\begin{proposition}\label{propositionForV-SIsolated}
If $S$ is a nontrivial  MPCS of tree graph $G$, then $\tilde{G}[\bar{S}]$ is a graph of isolated vertices.
\end{proposition}
\textbf{Proof}\quad Suppose that there exists a component $T$ of
$\tilde{G}[\bar{S}]$ that contains at least two vertices. See Fig.\ref{ProofForBarSIsolatedSet}.
%%%%%%%%%%%%%%  ?5 ?? %%%%%%%%%%%%
Then, $T$ is a tree graph and there exist at least two pendent vertices of $T$, say $u,v$.
Let $P_{uv}$ be the $(u,v)$-path in $T$ and $wv\in E(P_{uv})$. Considering the disconnected graph $\tilde{G}-wv$, it is seen that $u$ and $v$ belong to different components. Let $T_1$ be the component that contains $u$. Let $S_u=V(T_1)\cap S$ and set
\[S_1=S-S_u.\]
We claim that $S_1$ is a PCS. In fact, $S_u$ and $S_1$ belong to different components of $\tilde{G}-wv$. Hence, the submatrix $\textbf{L}_{S_1\rightarrow S_u}=\textbf{0}$. Further, set $\textbf{y}$ as the inducing eigenvector of MPCS $S$. Then, set $\tilde{\textbf{y}}=(\textbf{y}_{S_1},\textbf{0})$. For each $v'\in \bar{S}$, either $N_{\tilde{G}-wv}(v')\cap S\subseteq S_u$ or $N_{\tilde{G}-wv}(v')\cap S\subseteq S_1$.
 Therefore, $\tilde{\textbf{y}}$ is the inducing eigenvector of $S_1$.

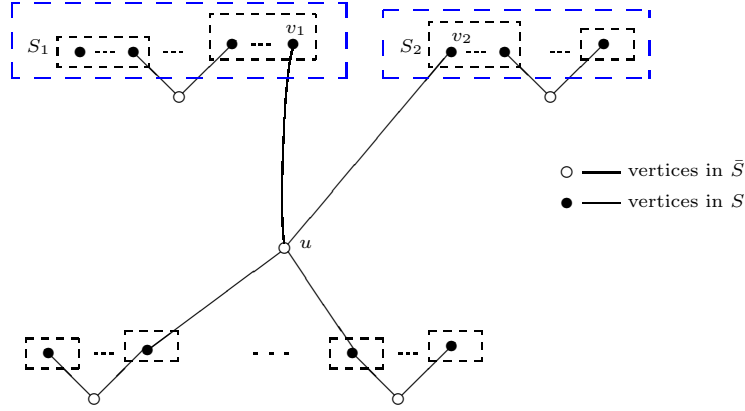
\begin{figure}[h]
\centering
\setlength{\unitlength}{1mm}
\begin{picture}(80,55)
\put(0,4){\dashbox{1}(7,4)}
\put(3,6){\circle*{1.5}}
\put(9,6){\line(1,0){0.5}}
\put(10,6){\line(1,0){0.5}}
\put(11,6){\line(1,0){0.5}}
\put(13,5){\dashbox{1}(7,4)}
\put(16,6.5){\circle*{1.5}}
\put(9,0){\circle{1.5}}
\put(9.5,0.5){\line(1,1){6}}
\put(8.5,0.5){\line(-1,1){6}}
\put(30,6){\line(1,0){0.5}}
\put(32,6){\line(1,0){0.5}}
\put(34,6){\line(1,0){0.5}}

\put(40,4){\dashbox{1}(7,4)}
\put(43,6){\circle*{1.5}}
\put(49,6){\line(1,0){0.5}}
\put(50,6){\line(1,0){0.5}}
\put(51,6){\line(1,0){0.5}}
\put(53,5){\dashbox{1}(7,4)}
\put(56,7){\circle*{1.5}}
\put(49,0){\circle{1.5}}
\put(49.5,0.5){\line(1,1){7}}
\put(48.5,0.5){\line(-1,1){6}}

\put(34,20){\circle{1.5}}
\put(36,20){\scriptsize{$u$}}

\put(33.5,19.5){\line(-4,-3){18}}
\put(34.5,19.5){\line(2,-3){9}}

\put(53,44){\dashbox{1}(12,6)}
\put(63,46){\circle*{1.5}}
\put(58,46){\line(1,0){0.5}}
\put(59,46){\line(1,0){0.5}}
\put(60,46){\line(1,0){0.5}}
\put(56,46){\circle*{1.5}}
\put(56,47.5){\scriptsize{$v_2$}}
\put(69,46){\line(1,0){0.5}}
\put(70,46){\line(1,0){0.5}}
\put(71,46){\line(1,0){0.5}}
\put(73,45){\dashbox{1}(7,4)}
\put(76,47){\circle*{1.5}}
\put(69,40){\circle{1.5}}
\put(69.5,40.5){\line(1,1){7}}
\put(68.5,40.5){\line(-1,1){6}}

\put(34.5,20.5){\line(5,6){21}}

\color{blue}{\put(47,42.5){\dashbox{2}(35,9)}}
\color{black}

\put(48,46){\scriptsize{$S_2$}}

\put(3,44){\dashbox{1}(12,4)}
\put(13,46){\circle*{1.5}}
\put(8,46){\line(1,0){0.5}}
\put(9,46){\line(1,0){0.5}}
\put(10,46){\line(1,0){0.5}}
\put(6,46){\circle*{1.5}}
\put(17,46){\line(1,0){0.5}}
\put(18,46){\line(1,0){0.5}}
\put(19,46){\line(1,0){0.5}}
\put(23,45){\dashbox{1}(14,6)}
\put(26,47){\circle*{1.5}}
\put(33,48.5){\scriptsize{$v_1$}}
\put(19,40){\circle{1.5}}
\put(34,47){\circle*{1.5}}
\put(28.5,47){\line(1,0){0.5}}
\put(29.5,47){\line(1,0){0.5}}
\put(30.5,47){\line(1,0){0.5}}
\put(19.5,40.5){\line(1,1){7}}
\put(18.5,40.5){\line(-1,1){6}}

\cbezier[500](32.8,20.6)(32,25)(33,49)(34.5,47)

\color{blue}
\put(-3,42.5){\dashbox{2}(44,10)}
\color{black}
\put(-1,46){\scriptsize{$S_1$}}

\put(70,30){\circle{1.5}}
\put(72,30){\line(1,0){5}}
\put(78,29.5){\scriptsize{vertices in $\bar{S}$}}

\put(70,26){\circle*{1.5}}
\put(72,26){\line(1,0){5}}
\put(78,25.5){\scriptsize{vertices in $S$}}

\end{picture}
\vspace{0.5cm}
\caption{Vertex $u\in\bar{S}$ with $d_{\tilde{G}}(u)>2$}
\label{ProofForVerticesDegreeInBarS}
\end{figure}

  Since $S_1$ is a PCS and $S_1\subset S$, it is contradictory to the fact that $S$ is an MPCS.\qed
\begin{proposition}\label{propositionForDegreeInBarS}
If $S$ is an MPCS of tree graph $G$, then $d_{\tilde{G}}(v)=2(\forall v\in \bar{S})$.
\end{proposition}
\textbf{Proof}\quad  From Lemma \ref{lemma1}, $d_{\tilde{G}}(v)\neq 1(\forall v\in \bar{S})$.
From Proposition \ref{propositionForV-SIsolated}, $N_{V(\tilde{G})}(v)\subset S(\forall v\in \bar{S})$.
Suppose that there exists a vertex in $\bar{S}$ , say $u$ , such that $d_{\tilde{G}}(u)\neq 2$, i.e. $d_{\tilde{G}}(u)> 2$. (see Fig.\ref{ProofForVerticesDegreeInBarS}).
Considering the subgraph $\tilde{G}-u$, different vertices that are adjacent to $u$ belong to different components of $\tilde{G}-u$. Let $T_i$ be the component and $v_i\in V(T_i)(i=1,2)$. Let $S_i=S\cap V(T_i)(i=1,2)$. Since $G$ is a tree graph, we have $S_1\cap S_2=\emptyset$. Let $\tilde{S}=S_1\cup S_2$. From $d_{\tilde{G}}(u)>2$, we know that $\tilde{S}$ is a proper subset of $S$.

We claim that $\tilde{S}$ is a PCS. In fact, let $\textbf{y}$ be the inducing eigenvector of MPCS $S$, where $y_i$ are its elements corresponding to $v_i(i=1,2)$ and $\textbf{y}_{S_i}$ are its elements corresponding to $S_i$. Let
\[\textbf{y}_{\tilde{S}}=\left[
\begin{array}{c}
\frac{\textbf{y}_{S_1}}{y_1}\\
\quad \\
-\frac{\textbf{y}_{S_2}}{y_2}
\end{array}
\right].
\]
Then, we have $\textbf{L}_{\tilde{S}
\rightarrow\tilde{S}}\textbf{y}_{\tilde{S}}=\lambda\textbf{y}_{\tilde{S}}$ and $\textbf{L}_{\bar{\tilde{S}}
\rightarrow\tilde{S}}\textbf{y}_{\tilde{S}}=\textbf{0}.$
Therefore, $[\textbf{y}_{\tilde{S}},\textbf{0}]$ is a inducing eigenvector for $\tilde{S}$ to be a PCS.
This is a contradiction to the fact that $S$ is an MPCS since $\tilde{S}$ is a proper subset of $S$.\qed

\subsection{No Isolated MPCS with $|S|\geq 3$ in Trees}
In this subsection, it is proved that there is no MPCS with isolated vertex set $S(|S|\geq 3)$ for tree graph.
\begin{theorem}\label{TheoremForIsolatedVertexSet}
Let $S$ be an isolated vertex set of tree graph $G$ and $|S|\geq 3$. Then, $S$ is not an MPCS.
\end{theorem}
\textbf{Proof}\quad Suppose that $S$ is an MPCS and it is an isolated vertex set. Then, $\textbf{L}_{S\rightarrow S}$ is a diagonal matrix. By $\textbf{L}_{S\rightarrow S}\textbf{y}_S=\lambda \textbf{y}_s$, we know that $d_G(v)=\lambda(\forall v\in S)$. Thus, all vertices in $S$ have the same degree.

On the other hand, from Proposition\ref{proposition4}, there exist at least two vertices in $S$ with degree 1.
Hence, we have
$d_G(v)=1 (\forall v\in S).$

Noting that $\tilde{G}[S\cup \bar{S}]$ is a tree graph and from Proposition \ref{propositionForDegreeInBarS}, we have
\[|S|+2|\bar{S}|=\sum d_i=2\epsilon=2(|S|+|\bar{S}|-1).\]
Hence, $|S|=2$. This is a contradiction to $|S|\geq 3.$
\qed

From Theorem \ref{TheoremForIsolatedVertexSet} and Theorem \ref{theorem3},
for the tree graph, the isolated vertex set $S$ is an MPCS if and only if $|S|=2$ and the two vertices in $S$ are twin nodes.

\subsection{A Type of Special MPCS of Trees}
Considering the Laplacian matrix of graph $\tilde{G}$, similar to the proof of Theorem \ref{TheoremForIsolatedVertexSet}, it is  seen that the following proposition holds.

\begin{proposition}\label{PropositonOfIsolatedVertexInG[S]}
 Let $S$ be a CS of tree graph $G$. If there exists a vertex $v_0\in S$ such that $d_{\tilde{G}[S]}(v_0)=0$, and $d_{\tilde{G}}(v_0)=t$, then the corresponding eigenvalue $\lambda=t$ and $d_{\tilde{G}}(v)=t(\forall v\in S: d_{\tilde{G}[S]}(v)=0 ).  \qed$.
\end{proposition}

\begin{proposition}\label{propositionOfPendentVertex}
 Let $S$ be an MPCS of tree graph $G$. If there exists a vertex $v_0\in S$ such that $d_{\tilde{G}[S]}(v_0)=0$, and $d_{\tilde{G}}(v_0)=1$ for each $ v $ with $d_{\tilde{G}}(v)=1 $ and $u_vv\in E(\tilde{G})$, then $u_v\in\bar{S}$.
\end{proposition}

\textbf{Proof}\quad By contrast, suppose that $u_v\in S$. Then, by Lemma \ref{lemma1}, we have $v\notin \bar{S}$. Hence, $v\in S$. From Proposition \ref{PropositonOfIsolatedVertexInG[S]}, we know that $\lambda=1$. Since $\textbf{L}_{v\rightarrow V(\tilde{G})}\textbf{y}=\lambda\textbf{y}_v=\textbf{y}_v$, $\textbf{y}_v-\textbf{y}_{u_v}=\textbf{y}_v$ and $\textbf{y}_{u_v}=0$. This is a contradiction to  $u_v\in S$ and $S$ is an MPCS. $\qed$

%8.S????????SS???????
To prove the next property of MPCS, we first give the following lemma.
\begin{lemma}\label{LemmaForNeighborInBarSLessThen2}
Let $\textbf{A}=(a_{ij})_{n \times n}$ be a Laplacian matrix of a tree graph. Matrix $\textbf{B}=(b_{ij})_{n \times n}$, where
$b_{ij}=\left\{
\begin{array}{ll}
a_{ij}+c_i,&  i=j,\\
a_{ij}, &  i\neq j,
\end{array}
\right. c_i\in\{-1,0\}$.
Then, there exists a vector $\textbf{X}_0=(x_1,x_2,\cdots,x_n)^\mathrm{T}(\forall i,x_i\neq 0)$ such that $\textbf{BX}_0=\textbf{0}$ if and only if $c_1=c_2=\cdots=c_n=0$.
\end{lemma}
\textbf{Proof}\quad $(\Leftarrow)$ If $c_1=c_2=\cdots=c_n=0$, then $\textbf{B}=\textbf{A}$ is a Laplacian matrix and $\textbf{X}_0=(1,1,\cdots,1)^{\mathrm{T}}$.
$(\Rightarrow)$ By induction on the order $n$.

If $n=1$, then $a_{11}=0$ and $(a_{11}+c_1)x_1=0$. Hence, $c_1=0$ since $x_1\neq 0$.

Assume that the lemma is true for each matrix of order lower than $n(n\geq 2)$.

Let $T$ be the  given tree with $n$ vertices and $A$ be its Laplacian matrix.  Let $v_1$ be a pendent vertex in $T$ and $v_2$ be the only  vertex adjacent to $v_1$. Then, matrix $A$ and matrix $B$ have the following form:
\[\textbf{A}=\left[
\begin{array}{ccccc}
1  &   -1  &  0 & \cdots  &  0\\
-1  & a_{22}  &  a_{23}   &  \cdots  &  a_{2n}\\
0  &  a_{32}  &  a_{33}  &  \cdots &  a_{3n}\\
\vdots  & \vdots  &  \vdots  &  \ddots  &  \vdots\\
0   &  a_{n2}  &   a_{n3}  & \cdots  & a_{nn}
\end{array}
\right],\]
\[\textbf{B}=\left[
\begin{array}{ccccc}
1+c_1  &   -1  &  0 & \cdots  &  0\\
-1  & a_{22}+c_2  &  a_{23}   &  \cdots  &  a_{2n}\\
0  &  a_{32}  &  a_{33}+c_3  &  \cdots &  a_{3n}\\
\vdots  & \vdots   &  \vdots   &  \ddots  &  \vdots \\
0   &  a_{n2}  &   a_{n3}  & \cdots  & a_{nn}+c_n
\end{array}
\right].\]

Let $\textbf{X}_0=(x_1,x_2,\cdots,x_n)(\forall i, x_i\neq 0)$ and $\textbf{BX}_0=\textbf{0}$.
From the first row of $\textbf{BX}_0=\textbf{0}$,
\[(1+c_1)x_1-x_2=0.\]
If $c_1=-1$, then $x_2=0$ and this is a contradiction. Hence, $c_1=0$ and $x_1=x_2$.

By elementary row transformation on matrix $\textbf{B}$,
\[\textbf{B}\xrightarrow{r_2:=r_1+r_2}
\left[
\begin{array}{cc}
1& \cdots \\
 \vdots  & \textbf{B}'
\end{array}\right].
\]
The second row of $\textbf{BX}_0=\textbf{0}$ and $x_1=x_2$ implies that
$\textbf{B}'\textbf{X}'_0=\textbf{0},$
where $\textbf{X}'_0=(x_2,x_3,\cdots,x_n)^{\mathrm{T}}.$
Noting that matrix $\textbf{B}'$ is of order $n-1$, we have $c_2=c_3=\cdots=c_n=0$.  $\qed$

$c_1=c_2=\cdots=c_n=0$ means $\mathbf{B}=\mathbf{A}$, that is $\mathbf{B}$ is the Laplacian matrix of tree graph.

Next, we will prove the following theorem to provide a graphical characterization of a special type of MPCS.
\begin{theorem}\label{theoremCharacterizationOfS}
Let $G$ be a tree graph and $S\subset V(G)$. If there exists a vertex $v_0\in S$ such that $d_{\tilde{G}[S]}(v_0)=0$, and $|N_{\tilde{G}}(v)\cap\bar{S}|\leq 1(\forall v\in S)$,
then $S$ is an MPCS if and only if $N_{\tilde{G}}(v)\cap \bar{S}=1(\forall v\in S)$.
\end{theorem}
\textbf{Proof}\quad $(\Rightarrow)$
We only need to prove that $d_{\tilde{G}}(v)=
d_{\tilde{G}[S]}(v)+1(\forall v\in S).$

 Since $S$ is an MPCS, $\tilde{G}$ is connected. Hence, for vertex $v_0$, from $|N_{\tilde{G}}(v)\cap \bar{S}|\leq 1$, we have $d_{\tilde{G}}(v_0)=1$.
By Proposition \ref{PropositonOfIsolatedVertexInG[S]}, the corresponding eigenvalue $\lambda=1$ and $d_{\tilde{G}}(v)=1(\forall v\in S: d_{\tilde{G}[S]}(v)=0 ).$ Therefor, only the vertices $v\in S$ with $d_{\tilde{G}[S]}(v)>0$ needed to be considered.

Let $v$ belong to $\tilde{G}_i$(tree), which is the $i$-th component of the subgraph $\tilde{G}[S]$. Let $S_i=V(\tilde{G}_i)$ and
let $\textbf{L}^{\tilde{G}}$ , $\textbf{L}^{\tilde{G}[S]}$ denote the Laplacian matrix of graph $\tilde{G}[S\cup \bar{S}]$, $\tilde{G}[S]$, respectively.

 Since $S$ is an MPCS, we know that there exists an eigenvector $\textbf{y}$ such that $\textbf{y}_v\neq 0(\forall v\in S)$ and
 $\textbf{L}^{\tilde{G}}_{S_i\rightarrow S_i}\cdot\textbf{y}_{S_i}=\lambda \textbf{y}_{S_i}=\textbf{y}_{S_i}(\lambda=1).$
Thus, $(\textbf{L}^{\tilde{G}}_{S_i\rightarrow S_i}-\textbf{I})\cdot \textbf{y}_{S_i}=\textbf{0}.$ By Lemma \ref{LemmaForNeighborInBarSLessThen2}, we have

\begin{equation}\label{eq3-2}
\textbf{L}^{\tilde{G}}_{S_i\rightarrow S_i}-\textbf{I}=\textbf{L}^{\tilde{G}[S]}_{S_i\rightarrow S_i},
\end{equation}
and this completes the proof of necessity.

$(\Leftarrow)$
If $v$ is a pendent vertex in $\tilde{G}$, then $d_{\tilde{G}[S]}(v)=0$. From Proposition \ref{proposition4} and Proposition \ref{PropositonOfIsolatedVertexInG[S]}, we know that the corresponding eigenvalue $\lambda=1$ for any PCS $S'(S'\subseteq S)$.

To prove the sufficiency, we need to prove that there exists only one linearly independent vector $\textbf{y}$ (corresponding eigenvalue $\lambda=1$) such that

\begin{equation}\label{eq4}
\textbf{L}^{\tilde{G}}\cdot
\textbf{y}=\textbf{y},
\end{equation}
\begin{equation}\label{eq5}
\textbf{y}_v\neq 0(\forall v\in S),
\end{equation}
\begin{equation}\label{eq6}
\textbf{y}_v=0(\forall v\in\bar{S}).
\end{equation}
We construct a new graph $H$ as follows. Let $u_i\in V(H)$ represent $\tilde{G}_i$, where $\tilde{G}_i$ is the $i$-th component of $\tilde{G}[S]$ and $w_j\in V(H)$, and let represent the $j$-th vertex in $\bar{S}$.
$E(H)=\{ u_iw_j | $ the $j$-th vertex in $\bar{S}$ is adjacent to one of vertices in $\tilde{G}_i\}$.
$H$ is a tree graph. Let $u_0$ correspond to $v_0$ and let it be a pendent vertex in $H$. Set \[\textbf{y}_{u_0}=1,\textbf{y}_v=0(\forall v \in \bar{S}).\]
According to Proposition \ref{propositionForG[S]Disconnected}, $u_i(i\geq 1)$ exists. $\forall u_i(i\geq 1)\in V(H)$, there exists only one path between $u_0$ and $u_i$. Let $l_i$ be the length of this path (i.e., the number of edges on this path). By Proposition \ref{propositionForDegreeInBarS}, $l_i$ is even. Set
\begin{equation}\label{eq7}
\textbf{y}_{v}=(-1)^{\frac{l_i}{2}}(\forall v \in V(\tilde{G}_i)).
\end{equation}
It is easily seen that the vector $\textbf{y}$ defined above satisfies  (\ref{eq5}) and (\ref{eq6}).

$\forall w_j\in V(H)$, by Proposition \ref{propositionForV-SIsolated} and Proposition \ref{propositionForDegreeInBarS}, we know that
there exist exactly two vertices, say $u_{i_1}$ and $u_{i_2}$, such that $u_{i_k}w_j\in E(H)(k=1,2)$. Hence, from (\ref{eq7}), we have
\[\textbf{L}^{\tilde{G}}_{v\rightarrow V(\tilde{G})}\cdot \textbf{y}=1+(-1)=0=\textbf{y}_{v}(\forall v\in\bar{S}).\]

$\forall u_i\in V(H)$, let $S_i=V(\tilde{G}_i)$. Since $\textbf{L}^{\tilde{G}[S]}_{S_i\rightarrow S_i}$ is the Laplacian matrix of tree graph $\tilde{G}_i$, we know that $\textbf{y}_{S_i}$(see (\ref{eq7})) is the only linearly independent vector such that
\[\textbf{L}^{\tilde{G}[S]}_{S_i\rightarrow S_i}\cdot\textbf{y}_{S_i}=\textbf{0}.\]
By (\ref{eq3-2}) and (\ref{eq1}), we have
\[\textbf{L}^{\tilde{G}}_{v\rightarrow V(\tilde{G})}\cdot\textbf{y}=\textbf{y}_{v}(\forall v\in S).\]
Hence, the vector $\textbf{y}$ defined above satisfies (\ref{eq4}). This completes the proof of sufficiency.  $\qed$

\section{Simulations and Discussion}\label{section4}

Finally, we illustrate the method for selecting leader vertices that render an undirected graph controllable by using the theories presented above.
\begin{figure}
\centering
\setlength{\unitlength}{1mm}
\begin{picture}(120,50)
\put(20,20){\circle*{1.5}}
\put(40,20){\circle{1.5}}
\put(60,20){\circle*{1.5}}
\put(80,20){\circle*{1.5}}
\put(100,20){\circle{1.5}}
\put(120,20){\circle*{1.5}}
\put(20,22){\scriptsize{$v_1$}}
\put(41,22){\scriptsize{$v_2$}}
\put(61,22){\scriptsize{$v_3$}}
\put(80,22){\scriptsize{$v_4$}}
\put(100,22){\scriptsize{$v_5$}}
\put(120,22){\scriptsize{$v_6$}}

\put(21,20){\line(1,0){18}}
\put(41,20){\line(1,0){58}}
\put(101,20){\line(1,0){18}}

\put(100,8){\circle{1.5}}
\put(102,8){\scriptsize{$v_7$}}
\put(100,9){\line(0,1){10}}

\put(40,30){\circle{1.5}}
\put(40,40){\circle{1.5}}
\put(40,21){\line(0,1){8}}
\put(40,31){\line(0,1){8}}
\put(41,30){\scriptsize{$v_8$}}
\put(41,40){\scriptsize{$v_{15}$}}

\put(60,30){\circle*{1.5}}
\put(60,40){\circle{1.5}}
\put(70,40){\circle*{1.5}}
\put(60,21){\line(0,1){18}}
\put(61,40){\line(1,0){8}}
\put(56,30){\scriptsize{$v_9$}}
\put(57,42){\scriptsize{$v_{13}$}}
\put(71,40){\scriptsize{$v_{14}$}}
\put(68,32){\scriptsize{$v_{10}$}}
\put(78,32){\scriptsize{$v_{11}$}}
\put(88,32){\scriptsize{$v_{12}$}}

\put(70,30){\circle*{1.5}}
\put(80,30){\circle{1.5}}
\put(90,30){\circle*{1.5}}
\put(61,30){\line(1,0){18}}
\put(81,30){\line(1,0){8}}

\end{picture}
\vspace{-0.8cm}
\caption{Tree graph $G$}
\label{Example}
\end{figure}
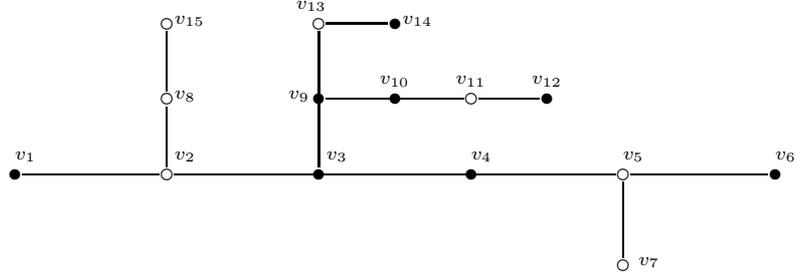

\subsection{ A Synthetic Tree}
See Fig.\ref{Example}. All MPCS of $G$ can be found, and the process of finding MPCS is shown in Table \ref{tableForExample}.

There are a total of three MPCS of $G$; they are
 $\{v_6,v_7\},
 \{v_6,v_4,v_3,v_1,$
 $v_9,v_{14},v_{10},v_{12}\},
   \{v_7,v_4,v_3,v_1,v_9,v_{14},v_{10},v_{12}\}.$
Hence, from Remark \ref{remarkForPropositionStatedWithMPVCS}, we know that the minimum number of leaders is 2 and there are a total of 15 different minimum leader  sets, i.e.,
\[\{v_6,v_7\}\,\, and \,\,\{v_i,v_k\}(i\in\{6,7\},k\in\{1,3,4,9,10,12,14\}).\]

\begin{table}
\centering
\caption{Process of finding MPCS}
\label{tableForExample}
\begin{tabular}{p{1.5cm}p{4cm}p{1.3cm}p{1.3cm}p{1.5cm}}
%{>{\columncolor{white}}lllll}
\toprule[1pt]
\rowcolor[gray]{0.94}
vertex  &  theory  &  $S$   &  $\bar{S}$  & MPCS\\
\midrule
\midrule
$v_6,v_7$  & Theorem \ref{theorem3} &  \checkmark  &  &\multicolumn{1}{>{\columncolor{white}[0pt][0pt]}l}{$\{v_6,v_7\}$}\\
\midrule
$v_6$ & initial(  or $v_7$)   & \checkmark &  &  \\

$v_7$ &  initial(or $v_6$)  &  & \checkmark &  \\

$v_5$ &  Lemma \ref{lemma1}(see $v_7$)  &  & \checkmark &  \\

$v_2 $ &  Proposition \ref{propositionOfPendentVertex}(see $v_5,v_6$)  &  & \checkmark &  \\

$v_8$ &  Proposition \ref{propositionOfPendentVertex}(see $v_5,v_6$)  &  & \checkmark &  \\

$v_{11}$ &  Proposition \ref{propositionOfPendentVertex}(see $v_5,v_6$)  &  & \checkmark &  \\

$v_{13}$ &  Proposition \ref{propositionOfPendentVertex}(see $v_5,v_6$)  &  & \checkmark &  \\

$v_4$ &  Proposition \ref{propositionForDegreeInBarS}(see $v_5$)  &  \checkmark & &  \\

$v_{15}$ &  Lemma \ref{lemma1}(see $v_8$)  &  & \checkmark &  \\

$v_3$ &  proof of Theorem \ref{TheoremForIsolatedVertexSet}(see $v_4,v_6$)  &  \checkmark &  &  \\

$v_{1}$ & Proposition \ref{propositionForDegreeInBarS}(see $v_2$)  & \checkmark & &  \\

$v_9$ &  Proposition \ref{propositionForDegreeInBarS},Theorem \ref{TheoremForIsolatedVertexSet}(see $v_{10}$)  & \checkmark & &  \\

$v_{14}$ & Proposition \ref{propositionForDegreeInBarS}(see $v_{13}$)  & \checkmark & &  \\

$v_{10}$ &  Lemma \ref{lemma1}(see $v_9,v_{11}$)  & \checkmark & &  \\

$v_{12}$ &  Proposition \ref{propositionForDegreeInBarS}(see $v_{11}$)  & \checkmark & &  \\

  &  Theorem \ref{theoremCharacterizationOfS} &   &  & $\{v_6\}\cup S'$ \\
   &   Theorem \ref{theoremCharacterizationOfS} &   &  & $\{v_7\}\cup S'$ \\

\bottomrule[1pt]
\end{tabular}
\flushleft \quad\,\, Note: $S'=\{v_4,v_3,v_1,v_9,v_{14},v_{10},v_{12}\}$
\end{table}

\subsection{DSFN}
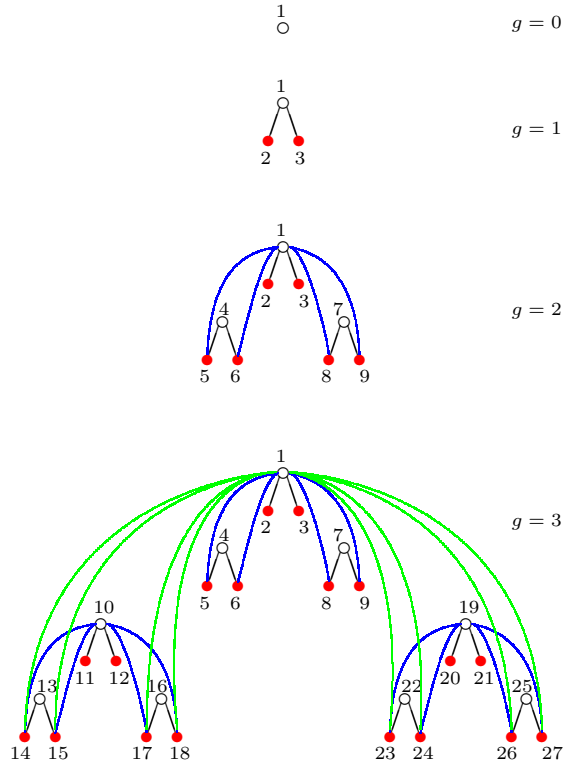
\begin{figure}
\centering
\setlength{\unitlength}{1mm}
\begin{picture}(120,100)
\put(60,100){\circle{1.5}}
\put(59,101.5){\scriptsize{$1$}}
\put(90,100){\scriptsize{$g=0$}}

\put(60,90){\circle{1.5}}
\put(59,91.5){\scriptsize{$1$}}
\put(90,86){\scriptsize{$g=1$}}
\put(58,85){\line(1,3){1.5}}
\put(62,85){\line(-1,3){1.5}}
\color{red}
\put(58,85){\circle*{1.5}}
\put(62,85){\circle*{1.5}}
\color{black}
\put(57,82){\scriptsize{$2$}}
\put(61.5,82){\scriptsize{$3$}}

%%%%%  g=2
\put(60,71){\circle{1.5}}
\put(59,72.5){\scriptsize{$1$}}
\put(90,62){\scriptsize{$g=2$}}
\put(58,66){\line(1,3){1.5}}
\put(62,66){\line(-1,3){1.5}}
\color{red}
\put(58,66){\circle*{1.5}}
\put(62,66){\circle*{1.5}}
\color{black}
\put(57,63){\scriptsize{$2$}}
\put(62,63){\scriptsize{$3$}}

\put(52,61){\circle{1.5}}
\put(51.5,62){\scriptsize{$4$}}
\put(50,56){\line(1,3){1.5}}
\put(54,56){\line(-1,3){1.5}}
\color{red}
\put(50,56){\circle*{1.5}}
\put(54,56){\circle*{1.5}}
\color{black}
\put(49,53){\scriptsize{$5$}}
\put(53,53){\scriptsize{$6$}}

\put(68,61){\circle{1.5}}
\put(66.5,62){\scriptsize{$7$}}
\put(66,56){\line(1,3){1.5}}
\put(70,56){\line(-1,3){1.5}}
\color{red}
\put(66,56){\circle*{1.5}}
\put(70,56){\circle*{1.5}}
\color{black}
\put(65,53){\scriptsize{$8$}}
\put(70,53){\scriptsize{$9$}}

\color{blue}
\cbezier[300](54,56.5)(55,60)(57,70)(59.3,71)
\cbezier[300](50,56.7)(50.5,70)(57,70.5)(59.3,71)
\cbezier[300](66,56.7)(65.5,60)(63,70.5)(60.8,71)
\cbezier[300](70,56.7)(69.5,70)(62,70.5)(60.8,71)
\color{black}

%%%%%%%%%%%%%%%%%
%%%%%%  g=3  %%%%%%%%%
\put(60,41){\circle{1.5}}
\put(59,42.5){\scriptsize{$1$}}
\put(90,34){\scriptsize{$g=3$}}
\put(58,36){\line(1,3){1.5}}
\put(62,36){\line(-1,3){1.5}}
\color{red}
\put(58,36){\circle*{1.5}}
\put(62,36){\circle*{1.5}}
\color{black}
\put(57,33){\scriptsize{$2$}}
\put(62,33){\scriptsize{$3$}}
\put(52,31){\circle{1.5}}
\put(51.5,32){\scriptsize{$4$}}
\put(50,26){\line(1,3){1.5}}
\put(54,26){\line(-1,3){1.5}}
\color{red}
\put(50,26){\circle*{1.5}}
\put(54,26){\circle*{1.5}}
\color{black}
\put(49,23){\scriptsize{$5$}}
\put(53,23){\scriptsize{$6$}}
\put(68,31){\circle{1.5}}
\put(66.5,32){\scriptsize{$7$}}
\put(66,26){\line(1,3){1.5}}
\put(70,26){\line(-1,3){1.5}}
\color{red}
\put(66,26){\circle*{1.5}}
\put(70,26){\circle*{1.5}}
\color{black}
\put(65,23){\scriptsize{$8$}}
\put(70,23){\scriptsize{$9$}}
\color{blue}
\cbezier[300](54,26.5)(55,30)(57,40)(59.3,41)
\cbezier[300](50,26.7)(50.5,40)(57,40.5)(59.3,41)
\cbezier[300](66,26.7)(65.5,30)(63,40.5)(60.8,41)
\cbezier[300](70,26.7)(69.5,40)(62,40.5)(60.8,41)
\color{black}
%%%%%%%
\put(36,21){\circle{1.5}}
\put(35,22.5){\scriptsize{$10$}}
\put(34,16){\line(1,3){1.5}}
\put(38,16){\line(-1,3){1.5}}
\color{red}
\put(34,16){\circle*{1.5}}
\put(38,16){\circle*{1.5}}
\color{black}
\put(32.5,13){\scriptsize{$11$}}
\put(37,13){\scriptsize{$12$}}
\put(28,11){\circle{1.5}}
\put(27.5,12){\scriptsize{$13$}}
\put(26,6){\line(1,3){1.5}}
\put(30,6){\line(-1,3){1.5}}
\color{red}
\put(26,6){\circle*{1.5}}
\put(30,6){\circle*{1.5}}
\color{black}
\put(24,3){\scriptsize{$14$}}
\put(29,3){\scriptsize{$15$}}
\put(44,11){\circle{1.5}}
\put(42,12){\scriptsize{$16$}}
\put(42,6){\line(1,3){1.5}}
\put(46,6){\line(-1,3){1.5}}
\color{red}
\put(42,6){\circle*{1.5}}
\put(46,6){\circle*{1.5}}
\color{black}
\put(40,3){\scriptsize{$17$}}
\put(45,3){\scriptsize{$18$}}
\color{blue}
\cbezier[300](30,6.5)(31,10)(33,20)(35.3,21)
\cbezier[300](26,6.7)(26.5,20)(33,20.5)(35.3,21)
\cbezier[300](42,6.7)(41.5,10)(39,20.5)(36.8,21)
\cbezier[300](46,6.7)(45.5,20)(38,20.5)(36.8,21)
\color{black}
%%%%%
\put(84,21){\circle{1.5}}
\put(83,22.5){\scriptsize{$19$}}
\put(82,16){\line(1,3){1.5}}
\put(86,16){\line(-1,3){1.5}}
\color{red}
\put(82,16){\circle*{1.5}}
\put(86,16){\circle*{1.5}}
\color{black}
\put(80.5,13){\scriptsize{$20$}}
\put(85,13){\scriptsize{$21$}}
\put(76,11){\circle{1.5}}
\put(75.5,12){\scriptsize{$22$}}
\put(74,6){\line(1,3){1.5}}
\put(78,6){\line(-1,3){1.5}}
\color{red}
\put(74,6){\circle*{1.5}}
\put(78,6){\circle*{1.5}}
\color{black}
\put(72,3){\scriptsize{$23$}}
\put(77,3){\scriptsize{$24$}}
\put(92,11){\circle{1.5}}
\put(90,12){\scriptsize{$25$}}
\put(90,6){\line(1,3){1.5}}
\put(94,6){\line(-1,3){1.5}}
\color{red}
\put(90,6){\circle*{1.5}}
\put(94,6){\circle*{1.5}}
\color{black}
\put(88,3){\scriptsize{$26$}}
\put(94,3){\scriptsize{$27$}}
\color{blue}
\cbezier[300](78,6.5)(79,10)(81,20)(83.3,21)
\cbezier[300](74,6.7)(74.5,20)(81,20.5)(83.3,21)
\cbezier[300](90,6.7)(89.5,10)(87,20.5)(84.8,21)
\cbezier[300](94,6.7)(93.5,20)(86,20.5)(84.8,21)
\color{black}
%%%%%%%%
\color{green}
\cbezier[500](26,7)(27,39)(58,41.1)(59,41.1)
\cbezier[500](30,7)(30,39)(58,41.1)(59,41.1)
\cbezier[500](42,7)(42,41)(58,41.1)(59,41.1)
\cbezier[500](46,7)(43,41)(58,41.1)(59,41.1)
\cbezier[500](74,7)(77,41)(60,41.1)(61,41.1)
\cbezier[500](78,7)(80,41)(60,41.1)(61,41.1)
\cbezier[500](90,7)(89,41)(59,41)(61,41.1)
\cbezier[500](94,7)(93,41)(59,41.1)(61,41.1)
\end{picture}
\caption{First four steps of construction of the DSFN
model $D(g)$}
\label{DSFN}
\end{figure}

A classical self-similar hierarchical network model called a deterministic scale-free network (DSFN)  was proposed in
\cite{Barabasi}. It is built in an iterative manner; see Fig.\ref{DSFN}.
Subsequently, in \cite{Kazumoto},
the authors studied the nature of the maximum eigenvalue in the
spectrum of DSFN and showed that most of the important quantities in network theory can be obtained analytically.
By using the graphical characterization of MPCS provided in this paper, we will determine both its minimum number and minimum sets of leaders in a more concise manner.

 Let $D(g)$, $n_l$, and $n_s$ denote the $g$-th generation of the DSFN , the minimum number of leaders, and the number of minimum sets of leaders, respectively.

By Theorem \ref{theorem3}, $\{3k-1,3k\}(k=1,2,\cdots,3^{g-1})$ are MPCS of $D(g)$. Except for these MPCS, $D(g)$ has no other MPCS. In fact,
see the vertices labeled $1,2,3$; if at least one of vertex $2,3$ belongs to $\bar{S}$, without loss of generality, let $2\in S$ and $3\in \bar{S}$. Then, $1\in \bar{S}$(Lemma \ref{lemma1}). Similarly, all vertices with label $3k-2(k=1,2,\cdots,3^g,g\geq 1)$ belong to $\bar{S}$.
Further, noting that all vertices labeled  $3k-1,3k$ are only adjacent to some vertices labeled  $3j-2$, one can find that $3k-1,3k(k=2,3\cdots,3^g,g\geq 1)$ belong to $\bar{S}$. Therefore, $S=\{2\}$, by (\ref{equationFor|S|>=2}), $S$ is not an MPCS.

These MPCS of $D(g)$ are pairwise disjoint; hence,
$n_l=3^{g-1}, n_s=2^{n_l}.$
The total number of vertices of $D(g)$ is $n=3^g$.

According to \cite{LiuAndSlotine}, the index $N_1$, defined as the ratio of $n_l$ to the network
size $n$, can be used to describe the level of network controllability.
For DSFN, the ratio is a constant,
$N_1=\frac{n_l}{n}=\frac{1}{3}.$
This indicates that it is necessary to
independently control at least one-third of the vertices to control the DSFN fully.

Since we have obtained all the MPCS of $D(g)$, we introduce another index $N_2$, which indicates the possibility of finding the minimum sets of leaders. It is defined as the ratio of $n_s$ to the number of subsets of $V$ with $n_l$ elements. For $D(g)$,
\begin{equation}\label{eqOfN2}
N_2(D(g))=\frac{n_s}{C_{n}^{n_l}}\leq(\frac{2}{3})^{n_l}\rightarrow 0(g\rightarrow \infty).
\end{equation}

Here, (\ref{eqOfN2}) suggests that it is extremely difficult to find the minimum set without theoretical analysis, even if we know the minimum number, and the probability rapidly decreases to 0.

\subsection{Cayley Tree}

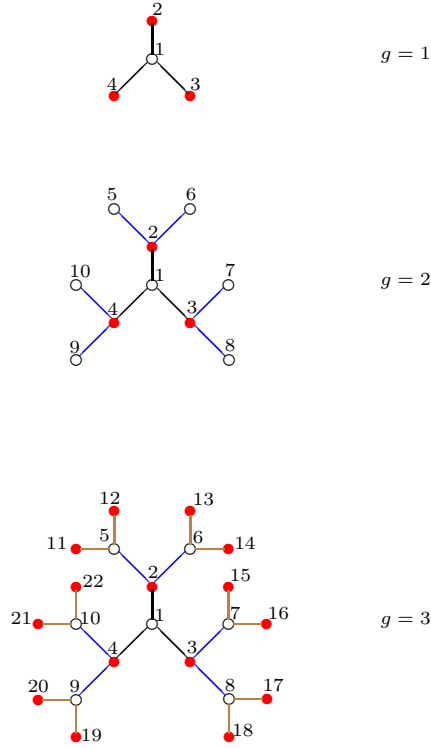
\begin{figure}
\centering
\setlength{\unitlength}{1mm}
\begin{picture}(120,100)
\color{red}
\put(60,100){\circle*{1.5}}
\put(55,90){\circle*{1.5}}
\put(65,90){\circle*{1.5}}
\color{black}
\put(60,100.8){\scriptsize{$2$}}
\put(54,90.9){\scriptsize{$4$}}
\put(65,90.8){\scriptsize{$3$}}
\put(90,95){\scriptsize{$g=1$}}
\put(60,95){\circle{1.5}}
\put(60.3,95.5){\scriptsize{$1$}}
\put(60,95.6){\line(0,1){4}}
\put(55.3,90.5){\line(1,1){4}}
\put(64.7,90.5){\line(-1,1){4}}
%%%%%%g=2
\put(60.3,65.5){\scriptsize{$1$}}
\put(59.4,71.2){\scriptsize{$2$}}
\put(64.5,61.2){\scriptsize{$3$}}
\put(54.1,61){\scriptsize{$4$}}
\put(54,76.2){\scriptsize{$5$}}
\put(64.4,76.2){\scriptsize{$6$}}
\put(69.5,66.2){\scriptsize{$7$}}
\put(69.5,56.2){\scriptsize{$8$}}
\put(49.1,56){\scriptsize{$9$}}
\put(49.1,66){\scriptsize{$10$}}

\color{red}
\put(60,70){\circle*{1.5}}
\put(55,60){\circle*{1.5}}
\put(65,60){\circle*{1.5}}
\color{black}
\put(90,65){\scriptsize{$g=2$}}
\put(60,65){\circle{1.5}}
\put(60,65.6){\line(0,1){4}}
\put(55.3,60.5){\line(1,1){4}}
\put(64.7,60.5){\line(-1,1){4}}
\color{black}
\put(55,75){\circle{1.5}}
\put(65,75){\circle{1.5}}
\color{blue}
\put(59.6,70.6){\line(-1,1){4}}
\put(60.3,70.6){\line(1,1){4}}
\color{black}

\put(50,65){\circle{1.5}}
\put(50,55){\circle{1.5}}
\color{blue}
\put(54.6,60.6){\line(-1,1){4}}
\put(50.4,55.6){\line(1,1){4}}
\color{black}

\put(70,65){\circle{1.5}}
\put(70.1,55.1){\circle{1.5}}
\color{blue}
\put(69.6,55.6){\line(-1,1){4}}
\put(65.4,60.6){\line(1,1){4}}
\color{black}
%%%%%%   g=3
\put(60.3,20.5){\scriptsize{$1$}}
\put(59.4,26.2){\scriptsize{$2$}}
\put(64.5,16.2){\scriptsize{$3$}}
\put(54.1,16){\scriptsize{$4$}}
\put(53,31){\scriptsize{$5$}}
\put(65.3,30.8){\scriptsize{$6$}}
\put(70.2,20.7){\scriptsize{$7$}}
\put(69.5,11.2){\scriptsize{$8$}}
\put(49.1,11){\scriptsize{$9$}}
\put(50.4,20){\scriptsize{$10$}}
\put(46,30){\scriptsize{$11$}}
\put(53,36){\scriptsize{$12$}}
\put(65.3,35.8){\scriptsize{$13$}}
\put(70.8,30){\scriptsize{$14$}}
\put(70.2,25.7){\scriptsize{$15$}}
\put(75.2,20.7){\scriptsize{$16$}}
\put(74.5,11.2){\scriptsize{$17$}}
\put(70.5,5.2){\scriptsize{$18$}}
\put(50.6,4.5){\scriptsize{$19$}}
\put(43.6,11){\scriptsize{$20$}}
\put(41.4,20){\scriptsize{$21$}}
\put(50.8,25){\scriptsize{$22$}}

\color{red}
\put(60,25){\circle*{1.5}}
\put(55,15){\circle*{1.5}}
\put(65,15){\circle*{1.5}}
\color{black}
\put(90,20){\scriptsize{$g=3$}}
\put(60,20){\circle{1.5}}
\put(60,20.6){\line(0,1){4}}
\put(55.3,15.5){\line(1,1){4}}
\put(64.7,15.5){\line(-1,1){4}}

\put(55,30){\circle{1.5}}
\put(65,30){\circle{1.5}}
\color{blue}
\put(59.6,25.6){\line(-1,1){4}}
\put(60.3,25.6){\line(1,1){4}}
\color{black}

\put(50,20){\circle{1.5}}
\put(50,10){\circle{1.5}}
\color{blue}
\put(54.6,15.6){\line(-1,1){4}}
\put(50.4,10.6){\line(1,1){4}}
\color{black}

\put(70,20){\circle{1.5}}
\put(70.1,10.1){\circle{1.5}}
\color{blue}
\put(69.6,10.6){\line(-1,1){4}}
\put(65.4,15.6){\line(1,1){4}}
\color{black}

%%%%
\color{red}
\put(50,30){\circle*{1.5}}
\put(55,35){\circle*{1.5}}
\color{brown}
\put(54.4,30){\line(-1,0){3.8}}
\put(55,30.6){\line(0,1){3.8}}
\color{black}

\color{red}
\put(70,30){\circle*{1.5}}
\put(65,35){\circle*{1.5}}
\color{brown}
\put(65,30.6){\line(0,1){3.8}}
\put(65.6,30){\line(1,0){3.8}}

\color{red}
\put(45,20){\circle*{1.5}}
\put(50,25){\circle*{1.5}}
\color{brown}
\put(50,20.6){\line(0,1){3.8}}
\put(49.4,20){\line(-1,0){3.8}}

\color{red}
\put(45,10){\circle*{1.5}}
\put(50,5){\circle*{1.5}}
\color{brown}
\put(50,9.4){\line(0,-1){3.8}}
\put(49.4,10){\line(-1,0){3.8}}

\color{red}
\put(70,25){\circle*{1.5}}
\put(75,20){\circle*{1.5}}
\color{brown}
\put(70.6,20){\line(1,0){3.8}}
\put(70,20.6){\line(0,1){3.8}}

\color{red}
\put(75.1,10.1){\circle*{1.5}}
\put(70.1,5.1){\circle*{1.5}}
\color{brown}
\put(70.7,10.1){\line(1,0){3.8}}
\put(70.1,9.3){\line(0,-1){3.8}}

\end{picture}
\caption{First three steps of construction of Cayley tree $C(g)$}
\label{CayleyTree}
\end{figure}

A Cayley tree is a classical model for dendrimers \cite{Julaiti}, which can be constructed in an iterative manner; see Fig.\ref{CayleyTree}. Let $C(g)$ denote the Cayley trees after $g$
generations. For any
$g\geq 2$, $C(g)$ is built from $C(g-1)$. For each pendent vertex
of $C(g-1)$, two new vertices are created and linked to the pendent.

\begin{figure}
\centering
\setlength{\unitlength}{1mm}
\begin{picture}(120,70)
\color{red}
\put(40,5){\circle{1.5}}
\put(45,5){\circle*{1.5}}
\put(50,5){\circle{1.5}}
\put(55,5){\circle*{1.5}}
\put(60,5){\circle{1.5}}
\put(65,5){\circle*{1.5}}
\put(70,5){\circle{1.5}}
\put(75,5){\circle*{1.5}}
\put(80,5){\circle{1.5}}
\put(85,5){\circle*{1.5}}
\put(90,5){\circle{1.5}}
\put(95,5){\circle*{1.5}}
\put(100,5){\circle{1.5}}
\put(105,5){\circle*{1.5}}
\put(110,5){\circle{1.5}}
\put(115,5){\circle*{1.5}}
\color{black}
\put(42.5,10){\circle{1.5}}
\put(52.5,10){\circle{1.5}}
\put(62.5,10){\circle{1.5}}
\put(72.5,10){\circle{1.5}}
\put(82.5,10){\circle{1.5}}
\put(92.5,10){\circle{1.5}}
\put(102.5,10){\circle{1.5}}
\put(112.5,10){\circle{1.5}}

\put(42,9.5){\line(-1,-2){1.8}}
\put(43,9.5){\line(1,-2){1.8}}
\put(52,9.5){\line(-1,-2){1.8}}
\put(53,9.5){\line(1,-2){1.8}}
\put(62,9.5){\line(-1,-2){1.8}}
\put(63,9.5){\line(1,-2){1.8}}
\put(72,9.5){\line(-1,-2){1.8}}
\put(73,9.5){\line(1,-2){1.8}}
\put(82,9.5){\line(-1,-2){1.8}}
\put(83,9.5){\line(1,-2){1.8}}
\put(92,9.5){\line(-1,-2){1.8}}
\put(93,9.5){\line(1,-2){1.8}}
\put(102,9.5){\line(-1,-2){1.8}}
\put(103,9.5){\line(1,-2){1.8}}
\put(112,9.5){\line(-1,-2){1.8}}
\put(113,9.5){\line(1,-2){1.8}}

\put(47.5,15){\circle*{1.5}}
\put(67.5,15){\circle*{1.5}}
\put(87.5,15){\circle*{1.5}}
\put(107.5,15){\circle*{1.5}}

\put(47.2,14.5){\line(-1,-1){4}}
\put(67.2,14.5){\line(-1,-1){4}}
\put(87.2,14.5){\line(-1,-1){4}}
\put(107.2,14.5){\line(-1,-1){4}}
\put(47.7,14.5){\line(1,-1){4}}
\put(67.7,14.5){\line(1,-1){4}}
\put(87.7,14.5){\line(1,-1){4}}
\put(107.7,14.5){\line(1,-1){4}}

\put(57.5,25){\circle*{1.5}}
\put(97.5,25){\circle*{1.5}}
\put(57.5,25){\line(-1,-1){10}}
\put(57.8,25){\line(1,-1){10}}
\put(97.5,25){\line(-1,-1){10}}
\put(97.8,25){\line(1,-1){10}}

\put(77.5,35){\circle*{1.5}}
\put(77.5,35){\line(-2,-1){20}}
\put(77.8,35){\line(2,-1){20}}

\put(47.5,50){\circle{1.5}}
\put(46.8,49.8){\line(-2,-1){29}}
\put(48.3,49.8){\line(2,-1){29}}
\color{red}
\put(17.5,35){\circle*{1.5}}
\put(47.5,60){\circle{1.5}}
\color{black}
\put(47.5,59.3){\line(0,-1){8.6}}
\put(47.5,60.6){\line(0,1){3}}
\put(17,34.5){\line(-2,-1){5}}
\put(17,34.5){\line(2,-1){5}}

\color{blue}
\put(38,0){\dashbox{2}(85,40)}
\color{green}
\put(5,0){\dashbox{2}(20,40)}
\put(34,55){\dashbox{2}(27,20)}
\color{black}

\put(99.5,2){\scriptsize{$4$}}
\put(104.5,2){\scriptsize{$3$}}
\put(109.4,2){\scriptsize{$2$}}
\put(114.5,2){\scriptsize{$1$}}
\put(103.7,9.1){\scriptsize{$5$}}
\put(113.5,10){\scriptsize{$6$}}
\put(108.4,15){\scriptsize{$7$}}
\put(85.2,15){\scriptsize{$8$}}
\put(98.4,25){\scriptsize{$9$}}
\put(53.8,25){\scriptsize{$10$}}
\put(78.8,35){\scriptsize{$11$}}
\put(48.8,50){\scriptsize{$12$}}
\put(13.4,34.8){\scriptsize{$13$}}
\put(48.8,60){\scriptsize{$14$}}

\put(12,20){\line(1,0){0.5}}
\put(15,20){\line(1,0){0.5}}
\put(18,20){\line(1,0){0.5}}

\put(45,70){\line(1,0){0.5}}
\put(48,70){\line(1,0){0.5}}
\put(51,70){\line(1,0){0.5}}

\end{picture}
\caption{Special MPCS of Cayley tree}
\label{findingMPCSsOfCayleyTree}
\end{figure}
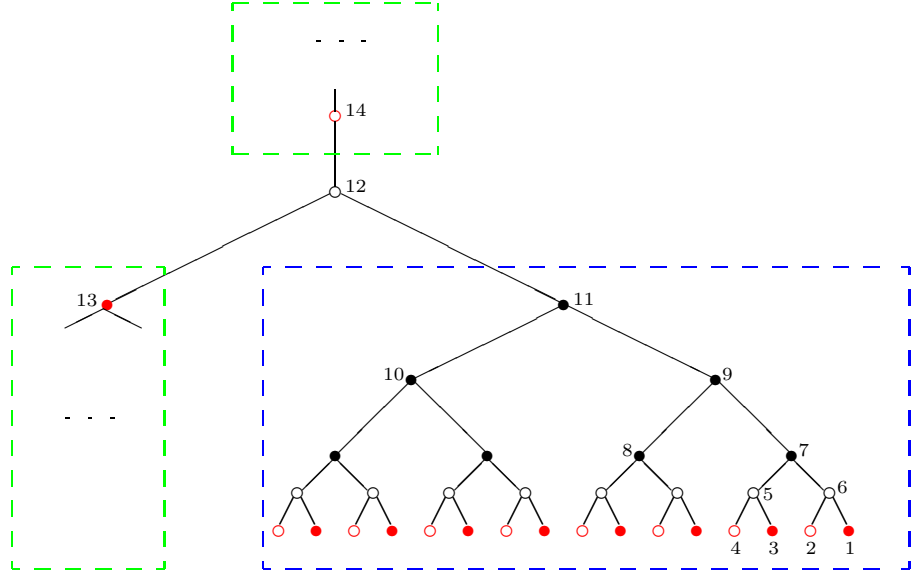

The total number of vertices of $C(g)$ is $n=3\cdot2^g-2$. Next, we will determine all MPCS of $C(g)$. Some MPCS can be determined by Theorem \ref{theorem3}. They are the sets composed of two pendents that are adjacent to the same vertex. In addition to these MPCS, there may be other MPCS. See Fig.\ref{findingMPCSsOfCayleyTree}.
Without loss of generality, let $S$ be an MPCS with $1\in S$ and $2\notin S$. Then, by Theorem \ref{theorem3}(see vertex 2), we have $6\notin S$. From Proposition \ref{propositionForDegreeInBarS}, $7\in S$. Further, $\lambda=1$ can be derived by Proposition \ref{PropositonOfIsolatedVertexInG[S]}. We claim that $9\in S$. Otherwise, the vertex $7$ will be an isolated vertex in $G[S]$ and recall the proof of Theorem \ref{TheoremForIsolatedVertexSet}; vertex 7 should have the same degree as vertex 1 in the Cayley tree. Similarly, $8\in S$. From $\textbf{L}_{7\rightarrow V}\textbf{y}_S=\lambda y_7$, we have $3y_7-y_9=y_7$, i.e., $y_9=2y_7$. Similarly, $y_9=2y_8$. From $\textbf{L}_{9\rightarrow V}\textbf{y}_S=\lambda y_9$, we have $y_9=y_{11}$, i.e., $11\in S$. Therefore, $10\in S$ and $y_{11}=y_{10}$. From $\textbf{L}_{11\rightarrow V}\textbf{y}_S=\lambda y_{11}$, we have $12 \notin S$.
Then, one of $\{13,14\}$ will belong to $S$. Thus far, we have obtained another MPCS $S$ of Cayley tree $C(g)(g\geq 5)$, where $S$ is composed of the vertices filled.

Table \ref{tableOfCayleyTree} lists  three important parameters $n_l, N_1,N_2$ of Cayley trees. From the last column of Table \ref{tableOfCayleyTree}, we can see that $N_2$ decrease to 0 rapidly, to the extent that the probability of finding the optimal solution is less than  one in ten million  after only 4 iterations.

\begin{table}
\centering
\caption{Minimum number of minimum sets of leaders of Cayley trees}
\label{tableOfCayleyTree}
\begin{tabular}
{>{\columncolor{white}}ccccc}
\toprule[1pt]
\rowcolor[gray]{0.94}
$g$ &  $n$  &   $n_l$  &  $N_1$   &  $N_2$\\
\midrule
\midrule
1  &  4  &  2  &  $\frac{1}{2}$    &   $\frac{1}{2}$ \\

\specialrule{0em}{4pt}{4pt}
2   &  10  &  3  &  $\frac{3}{10}$    &  $\frac{1}{15}$\\

\specialrule{0em}{4pt}{4pt}

3   &  22  &  6  &  $\frac{3}{11}$    &  $8.5776\!\!\times \!\!10^{-4}$\\

\specialrule{0em}{4pt}{4pt}

   4&  46  &  12  &  $\frac{6}{23}$    &  $1.0527\!\!\times \!\!10^{-7}$\\

\specialrule{0em}{4pt}{4pt}

  5&  94  &  26  &  $\frac{23}{47}$   & $2.2685\!\!\times\!\! 10^{-14}$   \\

\specialrule{0em}{4pt}{4pt}

  $g\geq 6 $  &   $3\!\cdot\!2^g\!\!-\!\!2$ &  $\frac{51}{64}\!\cdot\!2^g$   &  $\frac{17}{64}\,(g\rightarrow \infty)$   & $0\,\,(g\rightarrow \infty) $ \\

   $g\!\!\neq\!\! 5(mod 4)$  &   &    &    &   \\

\specialrule{0em}{4pt}{4pt}

$g\geq 6 $  &   $\quad 3\!\cdot\!2^g\!\!-\!\!2\quad $ &   $\quad\frac{27}{32}\!\cdot\!2^g\quad$   &  $\quad \frac{9}{32}\,(g\rightarrow \infty)$   & $0\,\,(g\rightarrow \infty) $ \\

   $g\!\!=\!\!5(mod 4)$  &   &    &    &   \\
\bottomrule
\end{tabular}
\end{table}

\section{Conclusion}\label{section5}
The controllability of undirected graphs has attracted considerable attention in recent years. However, the roles of leaders, the minimum number of leaders, and the leader location in a undirected graph are yet to be fully understood. A major contribution of this paper was to provide a method to determine the leaders directly from the topology structures of undirected graphs. These contributions can enhance our understanding of the leader's role in the controllability of undirected graphs.
To this end, we introduced the concepts of a critical set, perfect critical set, and minimal perfect critical set. These concepts indicate that some vertices with special graphical characterization should be selected as leaders. Necessary and sufficient conditions were proposed to uncover some special minimal perfect critical set. We described the graphical characterizations of minimal perfect critical 2  set. Further, we proved that minimal perfect critical 3 set does not exist. In addition, we pointed out that every tree graph is controllable when the leader set contains at least $|S_0|-1$ pendent vertices, where $S_0$ is the pendent vertex set.
Moreover, we studied an isolated
vertex set and found that if $S$ is isolated and $|S|\geq 3$, then $S$ is not an MPCS.
Finally, we completely described the MPCS of a tree graph with $|N_{\tilde{G}}(v)\cap\bar{S}|\leq 1(\forall v\in S)$ and showed that it can be used to discover some special MPCS of tree graphs.

All these results, which clearly indicate where the leaders are located, reveal the effect of the topology structure on controllability, and they could promote further investigation of the controllability of undirected graphs. To determine the characteristics of MPCS in a general network is an important topic that warrants in-depth analysis.

\end{document}